\documentclass[12pt]{article}
\usepackage{amsmath,amssymb}
\usepackage{amsthm}
\usepackage{amsfonts}   

\author{Simona Settepanella}
\title{\textbf{Cohomology of Pure Braid Groups of exceptional cases}}

\usepackage[all]{xy}
\usepackage{rotating}

\newcommand{\Q}{\mathbb Q}
\newcommand{\Z}{\mathbb Z}
\newcommand{\R}{\mathbb R}
\newcommand{\C}{\mathbb C}

\newcommand{\bW}{\mathbf W}

\newcommand{\bY}{\mathbf Y}

\newcommand{\aenne}{\mathbf{A_n}}
\newcommand{\benne}{\mathbf{B_n}}
\newcommand{\denne}{\mathbf{D_n}}

\newcommand{\fquattro}{\mathbf{F_4}}
\newcommand{\htre}{\mathbf{H_3}}
\newcommand{\hquattro}{\mathbf{H_4}}
\newcommand{\iduemm}{\mathbf{I_2}(m)}
\newcommand{\lgr}{\longrightarrow}

\newcommand{\A}{\mathcal A}

\newcommand{\sH}{\mathcal{H}_{(i_t,j_t)_{1\leq t\leq p}}}
\newcommand{\sG}{\mathcal{G}_{(i_t,j_t)_{1\leq t\leq p}}}

\begin{document}

\maketitle

{\small ABSTRACT. Consider the ring $R:=\Q[\tau,\tau^{-1}]$ of Laurent
  polynomials in the variable $\tau$. The Artin's Pure Braid Groups
  (or Generalized Pure Braid Groups) act over $R,$ where the action of
  every standard generator is the multiplication by $\tau$.  In this
  paper we consider the cohomology of these groups with coefficients
  in the module $R$ (it is well known that such cohomology is strictly
  related to the untwisted integral cohomology of the Milnor fibration
  naturally associated to the reflection arrangement). We compute this
  cohomology for the cases $\iduemm$, $\htre$, $\hquattro$,
  $\fquattro$ and $\aenne$ with $1 \leq \mathbf{n} \leq 7$.}

\section{Introduction}

Let $(\bW,S)$ be a finite Coxeter system realized as a reflection group 
in $\R^n$,
$\mathcal{A}(\bW)$ the 
arrangement in $\C^n$ obtained by complexifying the reflection hyperplanes
of $\bW$.Let
$$ 
\bY(\bW)=\bY(\mathcal{A}(\bW))=\C^n \setminus \cup_{H \in \mathcal{A}(\bW)}H.
$$ 
be the complement to the arrangement,
then $\bW$ acts freely on $\bY(\bW)$ and the fundamental group $G_W$ of 
the orbit space $\bY(\bW)/\bW$ is the so called \emph{Artin group} associated
to $\bW$ (see \cite{Bou}). Likewise the fundamental group $P_W$ of $\bY(\bW)$ 
is the \emph{Pure Artin group} or the pure braid group of the series $\bW$.
It is well known (\cite{Brie}) that these spaces $\bY(\bW)$ $(\bY(\bW)/\bW)$
are of type $K(\pi,1)$, so their cohomologies equal that of $P_W$ $(G_W)$.
   
The integer cohomology of $\bY(\bW)$ is well known (see \cite{Brie},\cite{OS},
\cite{b1},\cite{b5})
and so is the integer cohomology of the Artin groups associated to finite
Coxeter groups (see \cite{b9},\cite{b8},\cite{boss1}).

Let $R=\Q[\tau,\tau^{-1}]$ be the ring of rational Laurent polynomials.
To R can be given a structure of module over the Artin group $G_W$, where
standard generators of $G_W$ act as $\tau$-multiplication.
  
In \cite{b10} and \cite{b13} the authors compute the cohomology of all Artin 
groups associated to finite Coxeter groups with coefficients in the previous
module.

In a similar way we define a $P_W$-module $R_{\tau}$, where standard generators
of $P_W$ act over the ring $R$ as $\tau$-multiplication.

Equivalently, one defines an abelian local system (also called $R_{\tau}$) 
over $\bY(\bW)$ with fiber $R$ and local monodromy around each hyperplane
given by $\tau$-multiplication 
(for local systems on $\bY(\bW)$ see \cite{b14},\cite{b15}).

It is known that the cohomologies of $\bY(\bW)$ for the series 
$\aenne$, $\benne$ and $\denne$ stabilize, with respect to 
the natural inclusion, at a well known number of copies of the trivial 
$R$-module $\Q$ (\cite{Simona2}).
   
In this paper we give
that cohomology for the finite Coxeter groups  $\htre$, 
$\hquattro$, $\fquattro$ and $\aenne$ with 
$\mathbf{n}=7$ (for $\mathbf{n} \leq 5$ see \cite{CohenSuciu} and for $\mathbf{n}=6$ \cite{Denham2}).
We did these computations using an algorithm based on the Salvetti's complex.

Moreover, using results of G. Denham \cite{Denham} and H. Barcelo 
\cite{Barcelo}, we create 
algorithms to compute the cohomologies of the the flag complex 
$Fl_*(\bW)$ defined by Schechtman V. and Varchenko A. in \cite{SchVar} (see 
also \cite{Denham}). 
More precisely, these algorithms allow to perform calculations for finite Coxeter groups 
$\htre$, $\hquattro$, $\fquattro$ and $\aenne$ with $\mathbf{n} \leq 7$.

These computations are quite interesting and it's also very 
interesting to compare them. In fact this comparison
supports the conjecture that the integral cohomology of the Milnor Fibre 
associated to the reflection arrangements is torsion-free 
(see \cite{Denham} and \cite{Denham2}).

\section{Cohomology of pure braid groups $\htre$, $\hquattro$, $\fquattro$ and $\mathbf{A_7}$.}

In this section we give results obtained by direct computations.

Denote by $\varphi_i$ the cyclotomic polynomial having as roots the primitive
i-roots of $1$ and let 
$$
\{\varphi_i\}:=\Q[\tau,\tau^{-1}]/(\varphi_i)=\Q[\tau]/(\varphi_i)
$$
be the cyclotomic field of i-roots of $1$, thought as $R$-module.

\begin{table}[h]\label{tavola}
\begin{tabular}{|c|c|c|c|}
\hline
  & $\htre$ & $\hquattro$ & $\fquattro$ \\ 
\hline
  & & & \\
$H^0$ &  $0$ & $0$ & $0$ \\
\hline
  & & & \\
$H^1$ & $\{\varphi_1\}$   & $\{\varphi_1\}$  & $\{\varphi_1\}$ \\
\hline
  & & & \\
$H^2$ & $\{\varphi_1\}^{14}$   & $\{\varphi_1\}^{59}$   & $\{\varphi_1\}^{23}$  \\
\hline
  & & & \\
$H^3$ & $\{\varphi_1\}^{45} \oplus_{i \mid 15,i \neq 1} \{\varphi_i\}^{32}$ & $\{\varphi_1\}^{1079}\oplus \{\varphi_3\}$ & $\{\varphi_1\}^{167}\oplus \{\varphi_3\}^{8}$\\ 
\hline
  & & & \\
$H^4$ & $0$ & $\{\varphi_1\}^{6061}\oplus \{\varphi_3\}^{5039}$ & 
$\{\varphi_1\}^{385}\oplus 
\{\varphi_3\}^{232}$\\
  &  & $\oplus_{i \mid 60,i \neq 1,3 }\{\varphi_i\}^{5040}$ & 
$\oplus_{i \mid 24,i \neq 1,3}\{\varphi_i\}^{240}$\\
\hline
\end{tabular}
\caption{ Cohomologies of Pure Braid groups for exceptional cases.} 
\end{table}

\begin{table}
\begin{sideways}
\begin{tabular}{|r|r|r|r|r|r|r|r|}
\hline
  & $\mathbf{A_1}$ & $\mathbf{A_2}$ & $\mathbf{A_3}$ & $\mathbf{A_4}$ & $\mathbf{A_5}$ & $\mathbf{A_6}$ & $\mathbf{A_7}$ \\
\hline
  &   &  &  &  &  &  & \\
$H^0$ &  $0$ & $0$ & $0$ & $0$ & $0$ & $0$ & $0$\\
\hline
  &   &  &  &  &  &  & \\
$H^1$ & $\{\varphi_1\}$   & $\{\varphi_1\}$  & $\{\varphi_1\}$ & $\{\varphi_1\}$& $\{\varphi_1\}$ & $\{\varphi_1\}$ & $\{\varphi_1\}$\\
\hline
  &   &  &  &  &  &  & \\
$H^2$ & &$\{\varphi_1\}^{2} $ & $\{\varphi_1\}^{5}\oplus \{\varphi_3\}$   & $\{\varphi_1\}^{9}$ & $\{\varphi_1\}^{14}$ & $\{\varphi_1\}^{20}$ & $\{\varphi_1\}^{27}$ \\
  &   &$\oplus \{\varphi_3\}$  &  &  &  &  & \\
\hline
  &   &  &  &  &  &  & \\
$H^3$ &  &  & $\{\varphi_1\}^{6}\oplus \{\varphi_3\}$ 
&$\{\varphi_1\}^{26}\oplus \{\varphi_2\}^{2}$ & $\{\varphi_1\}^{71}\oplus\{\varphi_3\}$ & $\{\varphi_1\}^{155}\oplus\{\varphi_3\}$ &$\{\varphi_1\}^{295}$ \\
& & & $\oplus_{i\mid 6, i\neq 1}\{\varphi_i\}^2$ & & & & \\ 
\hline
  &   &  &  &  &  &  &  \\
$H^4$ &  & & &$\{\varphi_1\}^{24} \oplus\{\varphi_2\}^8$ &$\{\varphi_1\}^{154}\oplus\{\varphi_3\}^{14}$ & $\{\varphi_1\}^{580}$ & $\{\varphi_1\}^{1665}$\\
 & & & & $\oplus_{i\mid 10, i\neq 1,2}\{\varphi_i\}^6$ & $\oplus \{\varphi_5\}^6$ & $\oplus\{\varphi_3\}^{20}$ & \\
\hline
  &   &  &  &  &  &  & \\
$H^5$ & & & & &$\{\varphi_1\}^{120} \oplus \{\varphi_3\}^{37}$ & $\{\varphi_1\}^{1044}$ & $\{\varphi_1\}^{5104}$\\
 & & & & & $\oplus \{\varphi_5\}^30\oplus \{\varphi_{15}\}^{24}$& $\oplus\{\varphi_3\}^{121}$ &$ \oplus  \{\varphi_2\}^2$ \\
\hline
  &   &  &  &  &  &  & \\
$H^6$ & & & & & &$\{\varphi_1\}^{720}$ & $\{\varphi_1\}^{8028}$ \\
 & & & & & &$\oplus\{\varphi_3\}^{222}$ & $\oplus\{\varphi_2\}^{140}$\\
 & & & & & &$\oplus_{i\mid 21, i\neq 1,3}\{\varphi_i\}^{120}$ & $\oplus\{\varphi_7\}^{120}$\\
\hline
$H^7$ & & & & & & & $\{\varphi_1\}^{5040}$\\
 & & & & & & &$\oplus \{\varphi_2\}^{858} $ \\ 
 & & & & & & &$\oplus \{\varphi_7\}^{840}$ \\ 
 & & & & & & &$\oplus_{i \mid 28, i\neq 1,2,7} \{\varphi_i\}^{720}$ \\ 
\hline
\end{tabular}
\end{sideways}
\caption{ Cohomologies of Pure Braid Groups for $1 \leq \mathbf{n} \leq 7$.} 
\end{table}

\newpage

\section{Cohomology of $(Fl_*(\bW),\delta)$ for exceptional cases}

Recall, briefly, that a flag in a graded poset P is a tuple of elements 
$(X_0, \cdots , X_k)$ of P where each $X_i$ has rank $i$ and 
$X_i < X_{i+1}$ for $0 \leq i < k$ (the order given by the reverse inclusion). 
The flag complex $Fl_p$ in the intersection lattice of an arrangement $\A$, 
is defined to be the free abelian group of flags of length 
$p+1$, modulo the relations:
$$
\sum_{Y: X_{i-1} < Y < X_{i+1}} (X_0, \cdots , X_{i-1}, Y,  X_{i+1}, \cdots ,
X_p) = 0 ,
$$ 
where $i \geq 1$.

A set of $p$ independent hyperplanes $H_1,\cdots , H_p$ determines a flag with 
$X_i = H_1 \cap \cdots \cap H_i$ for $0 \leq i \leq p$ (see \cite{SchVar}). 
We will denote this 
flag by $\Phi=\lambda(H_1, \cdots ,H_p)$.
Accordingly, set $\varepsilon (H,\Phi) = (-1)^{i-1}$ and let define 
$\delta : Fl_p \lgr Fl_{p-1}$ as follows. 
Let $\Phi=(X_0, \cdots , X_p)$ be a flag and $H$ an hyperplane such that  
$H \leq X_p$. $D(\Phi,H)$ will denote the set of all $(p-1)$-flags 
\begin{equation}\label{insieme:Psi} 
\Psi = \lambda (H_1, \cdots , H_{i-1}, H_{i+1}, \cdots ,H_p),
\end{equation}  
such that $\Phi=\lambda(H_1,\cdots ,H_i=H,\cdots ,H_p)$.

Now, if $a : \A \lgr \Z_+$ is a weight function given by the multiplicities 
of the hyperplanes in an unreduced arrangement, we set (see \cite{SchVar})
\begin{equation} \label{delta}
\delta \Phi = \sum_{H \leq X_p} \sum_{\Psi \in D(\Phi,H)} 
\varepsilon (H, \Phi) a(H) \Psi .
\end{equation}

In the following table we give the cohomology of $(Fl_*(\bW),\delta)$
for reflection arrangements given by the finite Coxeter groups 
$\htre$, $\hquattro$ and $\fquattro$.

\begin{table}[h]\label{tavola2}
\begin{tabular}{|c|c|c|c|}
\hline
  & $\htre$ & $\hquattro$ & $\fquattro$ \\ 
\hline
  & & & \\
$H^0$ &  $0$ & $0$ & $0$ \\
\hline
  & & & \\
$H^1$ & $0$   & $0$  & $0$ \\
\hline
  & & & \\
$H^2$ & $0$   & $0$   & $0$  \\
\hline
  & & & \\
$H^3$ & $\oplus_{i \mid 15,i \neq 1} (\Z/i\Z)^{32}$ & $\Z/3\Z$ & $ (\Z/3\Z)^{8}$\\ 
\hline
  & & & \\
$H^4$ & $0$ & $(\Z/3\Z)^{5039}$ & $(\Z/3\Z)^{232}$\\
  &  & $\oplus_{i \mid 60,i \neq 1,3 }(\Z/i\Z)^{5040}$ & 
$\oplus_{i \mid 24,i \neq 1,3}(\Z/i\Z)^{240}$\\
\hline
\end{tabular}
\caption{ Cohomologies of flag complex for some exceptional cases.} 
\end{table}

\section{A filtration for the complex $(Fl_p(\aenne),\delta)$}

Recall that a NBC-basis for 
$A_p(\aenne)=A_p(\A(\aenne))$ is given by standard p-tuples  
$$
\{(H_{i_1,j_1},\cdots ,H_{i_p,j_p})\}_{1 \leq i_1 < \cdots <i_p \leq n}.
$$

and that
$\{\lambda(H_{i_1,j_1},\cdots ,H_{i_p,j_p})\}_{1 \leq i_1 < \cdots <i_p \leq n}$ is a basis for $Fl_p(\aenne)$ (see \cite{SchVar}).
 
In order to construct a filtration for the complex $(Fl_p(\aenne),\delta)$
we need some notations and definitions.

We consider the set $[n+1]:=\{1,\cdots ,n+1\}$, to each 
$H_{i,j}$ we associate the pair $(i,j)$;
then to each p-tuple $\sH:=(H_{i_1,j_1},\cdots ,H_{i_p,j_p})$ corresponds
 a graph 
$$
\sG:=G(\sH).
$$ 
For each point $h \in [n+1]$ of this graph we consider the set of vertices 
which are in the same connected component of $\sG$ as $h$; precisely: 
\begin{equation*}
\begin{split}
c_h(\sG):=\{i \in [n+1] \mid &\exists k_0=i,\cdots ,k_r=h \in [n+1] \\
&\textrm{ s.t. } H_{\sigma(k_q,k_{q+1})}\in \{\sH\}\textrm{ } \forall 
\textrm{ } 0 \leq q \leq r-1 \}
\end{split}
\end{equation*}
where
\begin{eqnarray}\label{mappa:sigma}
\sigma(i,j)= \left\{\begin{array}{ll}
 (i,j) & \textrm{if } i<j \\
 (j,i) & \textrm{if } i>j
 \end{array} \right. 
\end{eqnarray}
and 
$$
\{\sH\}:=\{H_{i_1,j_1},\cdots ,H_{i_p,j_p}\}.
$$

We set
$$
l_h(\sG):=\sharp c_h(\sG)
$$ 
its length.

These definitions can be extended to a flag 
$\Phi=\lambda(H_{i_1,j_1},\cdots ,H_{i_p,j_p})$. In this case we will 
denote:
\begin{equation*}
\begin{split}
&c_h(\Phi):=c_h(\sG) \\
&l_h(\Phi):=l_h(\sG)
\end{split}
\end{equation*}

Now let us consider the subcomplexes:
$$
G_{n+1}^k:=\{ \lambda(\sH) \in Fl_*(\aenne) \mid l_1(\sG) \leq k\},
$$
clearly the boundary map preserves $G_{n+1}^k$.

Let $F^k_{n+1}$ the cokernel of the natural inclusion in $Fl_*(\aenne)$
endowed
with the induced boundary.

Notice that the map $\delta$ defined on a flag 
$\Phi=\lambda(\sH)$ depends on the set $D(\Phi,H_{(i,j)})$ for 
$H_{(i,j)} \in \{\sH\}$ 
(see \ref{insieme:Psi}). In particular
it is easy 
to see that if $\Psi=\lambda(\mathcal{H}_{(i^{\prime}_t,j^{\prime}_t)_{1\leq t\leq p-1}}) \in D(\Phi,H_{(i,j)})$  
then for all $(i^{\prime},j^{\prime})$ s.t. $H_{(i^{\prime},j^{\prime})} \in \sH$ and $i^{\prime},j^{\prime} \notin c_i(\Phi)$, 
$H_{(i^{\prime},j^{\prime})} \in \{\mathcal{H}_{(i^{\prime}_t,j^{\prime}_t)_{1
\leq t\leq p-1}}\}$.

From the previous remark it follows that if  
$\Phi=\lambda(\mathcal{H}) \in (F^k_{n+1})_*$
is a flag with $l_1(\Phi)=k$ then exists a flag 
$\Phi^{\prime}=\lambda(\mathcal{H}^{\prime}) \in Fl_{*-(k-1)}(\aenne)$ s.t.
$$
\{\mathcal{H}^{\prime}\} \cup \{\mathcal{H}_{(i,j)}\}_{i,j \in c_1(\Phi)}=
\{\mathcal{H}\}
$$ 
and
$$
\delta_* \Phi=(\mathcal{H}_{(i,j)})_{i,j \in c_1(\Phi)}.\delta_{*-(k-1)}\Phi^{\prime}
$$
where the product is defined pointwise for all 
$\overline{H} \in \{\mathcal{H}_{(i,j)}\}_{i,j \in c_1(\Phi)}$ as follows:
$$
\overline{H}.\delta \Phi^{\prime}=\overline{H}.\sum_{H \leq X_p} 
\sum_{\Psi^{\prime} \in D(\Phi^{\prime},H)} 
\varepsilon (H, \Phi^{\prime}) a(H) \Psi^{\prime} :=\sum_{H \leq X_p} 
\sum_{\Psi^{\prime} \in D(\Phi^{\prime},H)} 
\varepsilon (H, \Phi^{\prime}) a(H) \overline{H}.\Psi^{\prime}
$$
and
$$
\overline{H}.\Psi^{\prime}=\overline{H}.\lambda(H_1,\cdots ,H_p):=
\lambda(H_1,\cdots H_{i-1},\overline{H}, H_i,\cdots ,H_p)
$$ 
with $H_{i-1} < \overline{H} < H_i$.

In other words for a flag $\Phi$ in $F^k_{n+1}$, $l_1(\Phi)=k$, the boundary 
map does not change components with indices in $c_1(\Phi)$.

Then, if we define a map of complexes
\begin{equation*}
\begin{split}
&i_n := i : Fl_*(\mathbf{A_{n-1}}) \lgr Fl_*(\aenne), \\
&i(\Phi)=i(\lambda(\sH))=\lambda(\mathcal{H}_{(i_t+1,j_t+1)_{1\leq t\leq p}})
\end{split} 
\end{equation*}

the cokernel of the map $i$ is the complex $F^2_{n+1}$ and we 
can iterate this construction considering the map
\begin{equation*}
\begin{split}
&i_n[1] := i : (n)^1 Fl_*(\mathbf{A_{n-2}})[1] \lgr F^2_{n+1}, \\
&i((j).\Phi)=i((j).\lambda(\sH))=
\lambda(\mathcal{H}_{(i^{\prime}_t,j^{\prime}_t)_{0\leq t\leq p}})
\end{split} 
\end{equation*}  
where $(n)^1=\{(i) \mid 2 \leq i \leq n \}$, 
$(i^{\prime}_0,j^{\prime}_0)=(1,j)$ and \\
$(i^{\prime}_t,j^{\prime}_t)=(2,\cdots,\check{j},\cdots,n+1)_{(i_t,j_t)}$, i.e.
the indexes in the positions $(i_t,j_t)$ in the $(n-1)$-tuple 
$(2,\cdots,\check{j},\cdots,n+1)$.

Notice that a chain of degree $k$ in $Fl_*(\mathbf{A_{n-2}})$ maps to one of 
degree 
$k+1$, for this we prefer to shift the complex 
$Fl_*(\mathbf{A_{n-2}})$ by $1$.

Again the cokernel of $i_n[1]$ 
is the complex $F^3_{n+1}$.

We continue in this way getting maps
\begin{equation*}
\begin{split}
&i_n[k-1] := i : (n)^{k-1} Fl_*(\mathbf{A_{n-k}})[k-1] \lgr F^{k}_{n+1},\\
&i((h_1,\cdots ,h_{k-1}).\Phi)=i((h_1,\cdots ,h_{k-1}).\sH)=\\
&(\mathcal{H}_{\sigma(h_i,h_{i+1})})_{0\leq i\leq k-1}.\lambda(\mathcal{H}_{(i^{\prime}_t,j^{\prime}_t)_{1 \leq t\leq p}})
\end{split} 
\end{equation*}  
where $h_0=1$, $\sigma$ is the map \ref{mappa:sigma},
$$
(n)^{k-1}=
\{(h_1,\cdots ,h_{k-1}) \mid 2 \leq h_i \leq n+1
\textrm{ and } h_i\neq h_j \textrm{ if } i\neq j\}
$$
and $(i^{\prime}_t,j^{\prime}_t)=(2,\cdots,\check{h_i},\cdots,n+1)_{(i_t,j_t)}$, i.e.
the indexes in the positions $(i_t,j_t)$ in the $(n-(k-1))$-tuple 
$(2,\cdots,\check{h_i},\cdots,n+1)$.

Each $i_n[k]$ gives rise to the exact sequence of complexes
\begin{equation} \label{eqn:suc.es.perflag}
0 \lgr (n)^{k} Fl_*(\mathbf{A_{n-k-1}})[k] \lgr
F^{k+1}_{n+1} \lgr F^{k+2}_{n+1}\lgr 0.
\end{equation}

The exact sequences \ref{eqn:suc.es.perflag} give rise to long exact sequences
in homology.

With this filtration we are able to compute the (co)-homology groups of $Fl_{\aenne}$ for $n \leq 7$ (see table \ref{tavola4}).

We can remark how, for these cases, the cohomology of the Flag Complex
is $0$ when the cohomology of the complement to the arrangement is a trivial
$\Z$-module. An interesting question could be to verify if this is always true.
In this case we could extend the ``stability'' theorem 
proved in \cite{Simona2} to the Cohomology of the Flag Complex.

\begin{table}\label{tavola4}
\begin{sideways}
\begin{tabular}{|r|r|r|r|r|r|r|}
\hline
  & $\mathbf{A_2}$ & $\mathbf{A_3}$ & $\mathbf{A_4}$ & $\mathbf{A_5}$ & $\mathbf{A_6}$ & $\mathbf{A_7}$ \\
\hline
     &  &  &  &  &  & \\
$H^0$  & $0$ & $0$ & $0$ & $0$ & $0$ & $0$\\
\hline
     &  &  &  &  &  & \\
$H^1$  & $0$ & $0$ & $0$& $0$ & $0$ & $0$\\
\hline
     &  &  &  &  &  & \\
$H^2$ & $\Z/3\Z$ & $\Z/3\Z$ & $0$& $0$ & $0$ & $0$\\
\hline
     &  &  &  &  &  & \\
$H^3$   &  &$\Z/3\Z$
&$(\Z/2\Z)^{2}$ & $\Z/3\Z$ & $\Z/3\Z$ & $0$ \\
 & &$\oplus_{i\mid 6, i\neq 1}(\Z/i\Z)^2$ & & & & \\
\hline
     &  &  &  &  &  &  \\
$H^4$   & & &$(\Z/2\Z)^8$ 
&$(\Z/3\Z)^{14}$& $(\Z/3\Z)^{30}$ & $0$ \\ 
  & & & $\oplus_{i\mid 10, i\neq 1,2}(\Z/i\Z)^6$ & $\oplus (\Z/5\Z)^6$  & & \\
\hline 
     &  &  &  &  &  & \\
$H^5$ & & & & $(\Z/3\Z)^{37}$ & & \\ 
& & & &$\oplus (\Z/5\Z)^{30}$ & 
$(\Z/3\Z)^{121}$ & $(\Z/2\Z)^{2}$\\
  & & & & $\oplus (\Z/15\Z)^{24}$ & & \\
\hline
     &  &  &  &  &  & \\
$H^6$  & & & & &$(\Z/3\Z)^{222}$ &
$(\Z/2\Z)^{140}\oplus(\Z/7\Z)^{120}$\\
  & & & & &$\oplus_{i\mid 21, i\neq 1,3}(\Z/i\Z)^{120} $ & \\
\hline
$H^7$  & & & & & &$(\Z/2\Z)^{858} \oplus (\Z/7\Z)^{840}$ \\
 & & & & & & $ \oplus_{i \mid 28, i\neq 1,2,7} (\Z/i\Z)^{720}$ \\ 
\hline
\end{tabular}
\end{sideways}
\end{table}

\phantom{\cite{b1,b5}}
\phantom{\cite{b8,b9,b14,b15}}
\phantom{\cite{Bou,boss1}}
\phantom{\cite{Brie,OS,Barcelo,Denham,SchVar}}
\phantom{\cite{b10,b13,CohenSuciu,Denham2,Simona2}}
\newpage{\pagestyle{empty}\cleardoublepage}
\bibliographystyle{plain}
\addcontentsline{toc}{chapter}{Bibliography}
\bibliography{maint} 

\end{document}